\documentclass{amsart}
\usepackage{amsmath,amscd,amssymb,graphicx}

\newtheorem{thm}{Theorem}[section]

\newtheorem{prop}[thm]{Proposition}
\newtheorem{lma}[thm]{Lemma}
\theoremstyle{definition}
\newtheorem{dfn}[thm]{Definition}

\def\N{\mathbb N}

\def\cP{\mathcal P}

\def\chra{\stackrel{c}{\hookrightarrow}}

\def\norm#1{\left\Vert#1\right\Vert}

\author[Tong]{ Simei Tong}
\address{
Simei Tong\\
Department of Mathematics, University of Wisconsin-Eau Claire\\Eau
Claire, WI 54701  } \email{tongs@uwec.edu} {\thanks{This paper was
a result in part during the 2009 summer Undergraduate Research
Experience in Pure and Applied Mathematics at the University of
Wisconsin - Eau Claire, supported by NSF-REU grant DMS-0552350 and
the Office of Research and Sponsored Programs of UW - Eau Claire
}}

\author[Phillipson]{Mitch Phillipson}
\address{
Mitch Phillipson\\
Department of Mathematics, University of Texas A \& M \\
College Station, TX 77843 } \email{phillippson@math.tamu.edu}

\author[DeFrain]{Isaac Defrain}
\address{
Isaac DeFrain\\
Florida Atlantic University Harriet L. Wilkes Honors College\\
Jupiter, FL 33458 } \email{idefrain@fau.edu}

\begin{document}
\title{Classifying Complemented Subspaces of $L_p$ with Alspach Norm}%

\dedicatory{This paper is dedicated to Professor Dale E. Alspach
on the occasion of his 60th Birthday. Professor Alspach was born in April 30,
1950. He was one of first three Ph.D. students to finish their
theses under the direction of William Johnson. Dr. Alspach's major
contribution to Banach Space Theory are the first example of a
non-expansive map on a weakly compact, closed convex subset of a
Banach space without a fixed point and classification results for
complemented subspaces of classical Banach spaces such as $C[0,1]$
and $L_p$, and translation invariant subspaces of $L_1(G)$ for $G$
Abelian. He is currently the Department Head at Oklahoma State
University. }

\begin{abstract}Understanding the complemented subspaces of $L_p$
has been an interesting topic of research in Banach space theory
since 1960. 1999, Alspach proposed a systematic approach to
classifying the subspaces of $L_p$ by introducing a norm given by
partitions and weights. This paper shows that with Alspach Norm we
are able to classify some complemented subspaces of $L_p,
2<p<\infty$.
\end{abstract}

\maketitle
\section{Introduction} Since 1960s understanding the complemented
subspaces of $L_p$ has been an interesting topic of research in
Banach space theory. Early in the work only obvious combinations
of $\ell_p$ and $\ell_2$ were known to give examples. In 1972,
Rosenthal's paper on sums of independent random variables was
seminal. He created $X_p$ and $B_p$ spaces in this paper. In 1975,
Schechtman proved that, up to isomorphism, there are infinitely
many complemented subspaces of $L_p$ by constructing tensor
products of $X_p$ spaces and in 1979, Bourgain, Rosenthal, and
Schechtman proved that, up to isomorphism, there are uncountably
many complemented subspaces of $L_p$. 1999, Alspach proposed a
systematic approach to classifying the subspaces of $L_p$ by
introducing a norm given by partitions and weights. His proposal
was the following:

Let A be a countable set and $P=\{ N_i\}$ be a partition of $A$
and
 $ W:A\rightarrow(0,1]$  be a function, which we refer to as
the \emph{weights}. Let $x_j\in\mathbb{R}$ for all $j\in A$.
Define
\begin{equation*}
\norm{(x_j)_{j\in A}}_{(P,W)} = \left(\sum_{N\in
P}\left(\sum_{j\in N} x_j^2
w_j^2\right)^\frac{p}{2}\right)^\frac{1}{p}
\end{equation*}
Suppose that $(P_k, W_k)_{k\in K}$ is a family of pairs of
partitions and weights as above. Define a (possibly infinite) norm
on the real-valued functions on $A$, $(x_i)_{i \in A},$ by
\begin{equation*}
 \norm{(x_i)} = \sup_{k\in K} \norm{(x_i)}_{(P_k,W_k)}
 \end{equation*}
and let $X$ be the subspace of elements of finite norm. In this
case we say that $X$ has an \emph{Alspach Norm}.\\

It is important to note that a space with the Alspach norm has a
natural unconditional basis and spaces with finite Alspach norm
are Banach spaces [AT1]. Alspach norm provides sequence space
realizations for some function spaces.  Alspach's approach gives a
unified description of many well known complemented subspaces of
$L_p$. It is proved that the class of spaces with such norms is
stable under $(p, 2)$ sums [AT1]. 2006, Alspach and Tong proved
that subspaces of $L_p$, $p>2$, with unconditional bases have
equivalent partition and weight norms[AT2]. In this article we
will explain what the conditions on partitions and weights will
produce certain known complemented subspaces of $L_p$. Our work is
far from complete due to the scope of this approach.
Classifying all complemented subspaces of $L_p$ with unconditional
basis is a big challenge and with Alspach
norm we made some progress.\\

From now on we will always assume that $p>2.$ In the rest of the
paper, we use Rosenthal's $X_p$ space and $B_p$ space many times,
[R]. Here are the definitions which can be found in Force's
dissertation [F]: $X_p$ can be realized as the closed linear span
in $L_p$ of a sequence $\{f_n\}$ of independent symmetric
three-valued random variables such that the ratios
$\norm{f_n}_2/\norm{f_n}_p$ approaches zero slowly. Another
realization of $X_p$ is as the set of all sequences $\{x_n\}$ in
$\ell_p$ for which the weighted $\ell_2$ norm $(\sum
|w_nx_n|^2)^{1/2}$ is finite and $(w_n)$ is a fixed sequence that
goes to zero slowly. The Banach space $B_p$ is of the form
$(Y_1\oplus Y_2 \oplus \cdots)_{\ell_p}$, where each space $Y_n$
is defined similar to $X_p$ but is isomorphic to $\ell_2$, and
$\{Y_n\}_{n=1}^{\infty}$ is chosen so that $\sup_{n\in \N}d(Y_n,
\ell_2)=\infty$, where $d(Y_n, \ell_2)$ is
the Banach-Mazur distance between $Y_n$ and $\ell_2$.\\

A partition is called \emph{discrete} if every subset has only one
element. A partition is called \emph{indiscrete} if the partition
is the whole set. A partition which is not discrete and not
indiscrete is called a \emph{regular} partition. Partition $P_1$
is called a \emph{refinement} of $P_2$ if every element in $P_2$
is a union of elements in $P_1$.  From now on, we treat the
discrete partition with constant weight of 1 as trivial and it
will be included in the discussion but will not be counted towards
the number of partitions. For instance we say that Rosenthal's
$X_p$ space [R] is an example of Banach space with Alspach norm given by
one partition and weights. In the context of subspaces of $L_p$
with unconditional basis there is always a lower $\ell_p$ estimate
and the discrete partition ensures that the spaces we consider
also have a lower $\ell_p$ estimate.

\section{One Regular Partition}

\begin{dfn}
 Let $A$ and $\cP=(P, W)$ be defined as above.
 A Banach space $X$ is said to have Alspach norm with one partition and weights if

\begin{equation*}
 \norm{(x_j)}_{X}=\max\left\{\left(\sum
x_j^p\right)^\frac1p,\left(\sum_{N\in P}\left(\sum_{j\in
N}x_j^2w_{j}^2\right)^\frac{p}2\right)^\frac1p\right\}
\end{equation*}
\end{dfn}
There is a complete classification of the spaces with Alspach norm
and one regular partition.

\begin{prop}
Let $P=\{N_i: i\in B\}$ where $B$ is an index set. Let $|N_i|$ be
the cardinality of $N_i$. Let $I=\{i: |N_i|=\infty\}$. Then
\begin{enumerate}
\item If $|I|< \infty$, then $X$ is isomorphic to one of $\ell_p$,
$X_p$, $\ell_2$, or $\ell_2\oplus \ell_p$.

\item If $|I|=\infty$, then $X$ is isomorphic to one of $\ell_p$,
$X_p$, $\ell_2\oplus \ell_p$, $B_p$, $(\sum \ell_2)_{\ell_p}$,
$(\sum \ell_2)_{\ell_p}\oplus X_p$, $B_p\oplus X_p$, or  $(\sum
X_p)_{\ell_p}$.
\end{enumerate}
\end{prop}

To prove the proposition we need following three lemmas. Below
$X_i=X_i^{|N_i|}=[e_j:j \in N_i]$ where $(e_j)$ is the natural
basis of $X$.

\begin{lma}
Let $I$ be defined as above. then
\begin{equation*}
X\sim \left (\sum_{i\in I}X_i\right )_{\ell_p} \oplus \left
(\sum_{i\notin I}X_i^{|N_i|}\right )_{\ell_p}
\end{equation*}
\end{lma}

\begin{lma}
If $B\diagdown I$ is finite, then
\begin{equation*}
\ell_p \oplus \left (\sum_{i\notin I}X_i^{|N_i|}\right
)_{\ell_p}\sim \ell_p
\end{equation*}
\end{lma}

\begin{lma}
If $B\diagdown I$ is infinite, then
\begin{equation*}
\left (\sum_{i\notin I}X_i^{|N_i|}\right )_{\ell_p}\sim \ell_p
\end{equation*}
\end{lma}

Proof: For each $i \notin I$, $X_i$ is a finite dimensional
version of one of the spaces considered by Rosenthal and thus is
isomorphic to a complemented subspace of $\ell_p$ and the norm of
the projection is independent of $i$. This implies
\begin{equation*}
\left (\sum X_i^{|N_i|}\right )_{\ell_p} \chra \left (\sum \ell_p
\right )_{\ell_p}.
\end{equation*}
Since $(\sum \ell_p)_{\ell_p}\sim \ell_p $, then
\begin{equation*}
\left (\sum X_i^{|N_i|}\right )_{\ell_p} \chra \ell_p.
\end{equation*}
Since every infinite dimensional complemented subspace of $\ell_p$
is isomorphic to $\ell_p$, then
\begin{equation*}
\left (\sum X_i^{|N_i|}\right )_{\ell_p} \sim {\ell}_p.
\end{equation*}

\medskip

After a messy computation based on splitting the argument into
several cases depending on the isomorphic type of the $\ell_p$ sum
of $X_i$ for $i\in I$ the results in the proposition follow.

\section{Admissible Partitions}

\begin{dfn}
A family of partitions and weights $\mathcal P$ is said to be
\emph{admissible} if there are partitions and weights $(P_0,W_0),$
$(P_1, W_1),$ and $(P_2,W_2)$ in $\mathcal P$ such that
$(P_0,W_0)$ is the discrete partition with weight constantly 1,
$(P_1, W_1)$ is a regular partition and weight and $P_2$ is the
indiscrete partition with weight $W_2=(w_{2,j})$.
\end{dfn}

We have following result

\begin{prop}
Assume $X$ be a sequence space of finite Alspach norm with an
admissible family of partitions and weights and only one regular
partition and weight. Then
\begin{enumerate}
\item If $\inf_j w_{2,j}\ge \delta >0$, then $X \thicksim \ell_2$.

\item Suppose $\sum_j (w_{2,j})^{\frac{2p}{p-2}} < {\infty}$. Let
$P_1=\{ N_i: i\in \N\}$. Let $|N_i|$ be the cardinality of $N_i$.
Let $I=\{i: |N_i|=\infty\}$. Then $X$ is isomorphic to one of the
spaces listed in Proposition 2.2 (1).

\item If we combine first two cases, i.e. there is some
$\delta>0$, such that $\{ j: w_{2,j} \ge \delta\}$ and $\{ j:
w_{2,j} < \delta\}$ are infinite and $\sum_{w_{2,j}<\delta}
(w_{2,j})^{\frac{2p}{p-2}} < {\infty}$, then $X$ is a direct sum
of $\ell_2$ from (1) and one of the spaces from (2).

%\item For every $\delta >0$, $\sum_{j:w_{2,j}<\delta}
%w_{2,j}^{\frac{2p}{p-2}}=\infty$, then $X$ is isomorphic to ???
\end{enumerate}
\end{prop}

Proof:
\begin{enumerate}
\item Since $\inf w_{2j}\ge \delta >0$, then
\begin{equation*}
\delta \left ( \sum_j x_j^2\right )^{\frac12} \le \norm {(x_j)}_X
\le \left ( \sum_j x_j^2\right )^{\frac12}
\end{equation*}
\item Since  $\sum_j (w_{2,j})^{\frac{2p}{p-2}} < {\infty}$, then
we can apply H\"older's inequality
\begin{equation*}
\left ( \sum x_j^2 w_{2,j}^2\right )^{\frac 12}\le \left (\sum
x_j^p \right )^{\frac 1p}\left (\sum
w_{2,j}^{\frac{2p}{p-2}}\right )^{\frac{p-2}{p}}
\end{equation*}
Thus
\begin{equation*}
\norm { (x_j)}_X \sim  \norm { (x_j)}_{P_1}
\end{equation*}
Apply the results from proposition 2.2 to the rest of the proof of
(2).
\end{enumerate}

\section{Two Partitions Related by Refinement} Let $A$ be any index set.
Let $P_1$ and $P_2$ be two partitions of $A$
and let $W_1=(w_{1,k})$ and $W_2=(w_{2,k})$ be two sequences of
weights.  Assume that $P_1$ is a
refinement of $P_2$ and that $w_{1,j}\ge w_{2,j}$ for all
$j\in A$. For a fixed $N \in P_1$, notice that
\begin{align*}
W_N &= \sup\left \{\left (\sum x_j^2 w_{2,j}^2\right
)^{\frac{1}{2}}/\left (\sum x_j^2 w_{1,j}^2\right )^{\frac{1}{2}}\right \} \\
&\le \sup_{k \in N} \frac {w_{2,k}}{w_{1,k}}
\sup \left \{\left (\sum x_j^2 w_{1,j}^2 \right )^{\frac{1}{2}}/\left (\sum x_j^2 w_{1,j}^2\right )^{\frac{1}{2}} \right \}\\
&= \sup_{k \in N} \frac{w_{2,k}}{w_{1,k}}
\end{align*}

the  supremum in the first two lines are taken over all sequences
$(x_j)$ for $x_j=0, j \notin N $ and $x_j \ne 0 $ for finitely
many $j$. Taking $x_j =1$ and $x_k=0$ for $k \ne j$,
shows that $\begin{displaystyle}W_N =
\sup_{k \in N}
\frac{w_{2,k}}{w_{1,k}}\end{displaystyle}$.\\

Now notice that for fixed $M \in P_2$,
\begin{align*}
\sum_{N \subset M} \sum_{j \in N} x_j^2 w_{2,j}^2 &\le \sum_{N
\subset M} W_N^2 \sum_{j \in N} x_j^2 w_{1,j}^2 \\ &\le \left
(\sum_{N \subset M} W_N^{\frac {2p}{p-2}}\right )^{\frac{p-2}{p}}
\left(\sum_{N \subset M} \left( \sum_{j \in N}
x_j^2 w_{1,j}^2\right )^{p/2}\right )^{2/p}
\end{align*}

Then by Using H\"older's inequality
\begin{align*}
\norm{(x_j)}&=\max\left\{\left(\sum
x_j^p\right)^\frac1p,\left(\sum_{N\in P_1}\left(\sum_{i\in
N}x_i^2w_{1i}^2\right)^\frac{p}2\right)^\frac1p , \left(\sum_{M\in
P_2}\left(\sum_{j\in
M}x_j^2w_{2j}^2\right)^\frac{p}2\right)^\frac1p\right\}\\
&=\max\left\{\left(\sum x_j^p\right)^\frac1p,\left(\sum_{N\in
P_1}\left(\sum_{i\in
N}x_i^2w_{1i}^2\right)^\frac{p}2\right)^\frac1p , \left(\sum_{M\in
P_2}\left(\sum_{N\subset M }\sum_{j\in
N}x_j^2w_{2j}^2\right)^\frac{p}2\right)^\frac1p\right\} \\
&\le \max\left\{\left(\sum x_j^p\right)^\frac1p,\left(\sum_{N\in
P_1}\left(\sum_{i\in
N}x_i^2w_{1i}^2\right)^\frac{p}2\right)^\frac1p , \left(\sum_{M\in
P_2}\left(\sum_{N\subset M }W_N^2\sum_{j\in
N}x_j^2w_{1j}^2\right)^\frac{p}2\right)^\frac1p\right\}\\
&\le \max\left\{\left(\sum x_j^p\right)^\frac1p,\left(\sum_{N\in
P_1}\left(\sum_{i\in
N}x_i^2w_{1i}^2\right)^\frac{p}2\right)^\frac1p , \left(\sum_{M\in
P_2}\left(\sum_{N\subset M }W_N^{\frac{2p}{p-2}}\right
)^{\frac{p-2}{2}}\left (\sum_{N \subset M}\sum_{j\in
N}x_j^2w_{1j}^2\right)^\frac{p}2\right)^\frac1p\right\}
\end{align*}

We have following result:

\begin{prop}
Let partition $P_1$ be a refinement of partition $P_2$. Let
$W_1=(w_{1,k})$ and $W_2=(w_{2,k})$ be two corresponding sequences
of weights. Assume that $w_{1,j}\ge w_{2,j}$ for all $j\in A$

If
\begin{equation*}
\sup_M\left\{\sum_{N\subset M }W_N^{\frac{2p}{p-2}}\right \}
<\infty
\end{equation*}
then $X$ can be classified by the behavior of $W_1$.

\end{prop}

\begin{prop} Let $X$ be the Banach space with finite norms given by two
partitions $P_1$ and $P_2$ and we also assume for every $N\in P_1$
there is an $M\in P_2$ such that $N\subset M$ and $w_{1k}\ge
w_{2k}$ for all $k\in N$. We have following results:

\begin{enumerate}
\item Assume $\inf(w_{2k})=\delta>0$
\begin{enumerate}
\item If $|P_2|<\infty$ we have
\begin{equation*}
  \delta|P_2|^{\frac1p-\frac12}\norm{x_j}_2\le \delta \left(\sum_{M\in
P_2} \left(\sum_{j\in M} x_j^2 \right)^\frac{p}2 \right)^\frac1p
\le \norm{(x_j)}_X\le\norm{(x_j)}_2.
\end{equation*}
 Therefore,
\begin{equation*}
 X\sim\ell_2.
\end{equation*}

\item If $|P_2|=\infty$ and $|M|<\infty$ for all $M\in P_2$,
\begin{equation*}
 \norm{x_j}_X\sim \left(\sum_{M \in P_2}\left(\sum_{j\in
M}x_j^2\right)^\frac{p}2\right)^\frac1p.
\end{equation*}
Let $Y_n$ represent the space generated by the inner sum which
implies that $X\sim (\sum Y_n)_p$.  Since $X$ is infinite
dimensional we conclude
\begin{equation*}
 X\sim\ell_p.
\end{equation*}
\item If $|P_2|=\infty$, $|M|=\infty$ for at least one and at most
finitely many $M\in
P_2$ Then $X$ will be a direct sum of the spaces,
\begin{equation*}
 X\sim\ell_2\oplus\ell_p.
\end{equation*}
\item If $|P_2|=\infty$ and $|M|=\infty$ for infinitely many $M\in
P_2$, then
\begin{equation*}
 \norm{x_j}_X\sim \left(\sum_{|M|=\infty}\left(\sum_{j\in
M}x_j^2\right)^\frac{p}{2}+\sum_{|M|<\infty}\left(\sum_{j\in
M}x_j^2\right)^\frac{p}2\right)^\frac1p.
\end{equation*}
and \begin{equation*}
 X\sim\left(\sum\ell_2\right)_{\ell_p}.
\end{equation*}

\end{enumerate}

\item Assume $\inf(w_{1k})=\gamma>0$ and $\inf(w_{2k})=0$
\begin{enumerate}
\item If $|P_1|<\infty$, we have

\begin{equation*}
\gamma|P_1|^{\frac{1}{p}-\frac{1}{2}}\norm{x_j}_2\le \gamma
\left(\sum_{N\in P_1} \left(\sum_{j\in N} x_j^2 \right)^\frac{p}2
\right)^\frac{1}{p} \le \norm{(x_j)}_X\le\norm{(x_j)}_2
\end{equation*}

so

\begin{equation*} X \sim \ell_2.
\end{equation*}

\item If $|P_1|= \infty$ and $|N|< \infty$ for all $N \in P_1$,
and for every $M\in P_2$, let $C_M= \sum_{n\in
M}w_{2n}^\frac{2p}{p-2}<\infty$ and assume $\sup_{M\in
P_2}C_M<\infty$, then

\begin{equation*}
 X\sim\ell_p
\end{equation*}

\end{enumerate}

\end{enumerate}
\end{prop}

\end{document}